# Networked Microgrids for Improving Economics and Resiliency

Guodong Liu, Thomas B. Ollis, Bailu Xiao
Power & Energy Systems Group
Oak Ridge National Laboratory
Oak Ridge, USA
Email: liug@ornl.gov, ollistb@ornl.gov, xiaob@ornl.gov

Xiaohu Zhang, Kevin Tomsovic
Dept. of Electrical Engineering and Computer Science
The University of Tennessee
Knoxville, USA
Email: xzhang46@utk.edu, tomsovic@utk.edu

***Abstract*—In this paper, we propose networked microgrids to facilitate the integration of variable renewable generation and improve the economics and resiliency of electricity supply in microgrids. A new concept, probability of successful islanding (PSI) is used to quantify the islanding capability of a microgrid considering the uncertainty of renewable energy resources and load as well as exchanged power at PCC. With the goal of minimizing the total operating cost while preserving user specified PSI, a chance-constrained optimization problem is formulated for the optimal scheduling of both individual microgrid and networked microgrids. Numerical simulation results show significant saving in electricity cost can be achieved by proposed networked microgrids without compromising the resiliency. The impact of correlation coefficients among the renewable generation and load of adjacent microgrids has been studied as well.**

***Index Terms*—Networked microgrids, optimal scheduling, probability of successful islanding, economics, resiliency.**

## NOMENCLATURE

The main symbols used in this paper are defined below. Others will be defined as required in the text. A $\triangle$ indicates forecast error for the variable while ˆ indicates the forecast value.

### A. Indices

| | |
|---|---|
| $n$ | Index of microgrids, running from 1 to $N_M$. |
| $i$ | Index of dispatchable generators, running from 1 to $N_G$. |
| $j$ | Index of demands, running from 1 to $N_D$. |
| $b$ | Index of battery storage devices, running from 1 to $N_B$. |
| $t$ | Index of time periods, running from 1 to $N_T$. |
| $m$ | Index of energy blocks offered by generators, running from 1 to $N_I$. |
| $l$ | Index of probability intervals, running from 1 to $N_L$. |

### B. Variables

*1) Binary Variables:*

| | |
|---|---|
| $u_{it}$ | 1 if unit $i$ is scheduled on during period $t$ and 0 otherwise. |
| $u_{bt}^{\mathrm{C}}$, $u_{bt}^{\mathrm{D}}$ | 1 if battery $b$ is scheduled charging/discharging during period $t$ and 0 otherwise. |
| $b_{tl}^{\mathrm{U}}$, $b_{tl}^{\mathrm{D}}$ | Binary indicators of probability interval $l$ during period $t$ . |

*2) Continuous Variables:*

| | |
|---|---|
| $p_{it}\left(m\right)$ | Power output scheduled from the $m$-th block of energy offer by dispatchable unit $i$ during period $t$ . Limited to $p_{it}^{\max}\left(m\right)$. |
| $P_{it}$ | Power output scheduled from dispatchable unit $i$ during period $t$. |
| $P_t^{\mathrm{PCC}}$ | Exchanged power at PCC during period $t$. |
| $P_{bt}^{\mathrm{C}}, P_{bt}^{\mathrm{D}}$ | Charging/discharging power of battery $b$ during period $t$. |
| $P_{bt}$ | Output power of battery $b$ during period $t$. |
| $SOC_{bt}$ | State of charge of battery $b$ during period $t$. |
| $R_{it}^{\mathrm{U}}$, $R_{it}^{\mathrm{D}}$ | Up- and down-spinning reserve of unit $i$ during period $t$. |
| $R_{bt}^{\mathrm{U}}$, $R_{bt}^{\mathrm{D}}$ | Up- and down-spinning reserve of battery $b$ during period $t$. |
| $\mathrm{PSI}_t$ | Probability of successful islanding during period $t$. |

### C. Constants

| | |
|---|---|
| $\lambda_{it}\left(m\right)$ | Marginal cost of the $m$-th block of energy offer by dispatchable unit $i$ during period $t$. |
| $C_{bt}$ | Degradation cost of battery $b$ during period $t$. |
| $\lambda_t^{\mathrm{PCC}}$ | Purchasing/selling price of energy from/to distribution grid during period $t$. |
| $A_i$ | Operating cost of dispatchable unit $i$ at the point of $P_i^{\min}$. |
| $Q_{it}^{\mathrm{U}}$, $Q_{it}^{\mathrm{D}}$ | Cost of up- and down-spinning reserve of unit $i$ during period $t$. |

This manuscript has been authored by UT-Battelle, LLC under Contract No. DE-AC05-00OR22725 with the U.S. Department of Energy. The United States Government retains and the publisher, by accepting the article for publication, acknowledges that the United States Government retains a non-exclusive, paid-up, irrevocable, world-wide license to publish or reproduce the published form of this manuscript, or allow others to do so, for United States Government purposes. The Department of Energy will provide public access to these results of federally sponsored research in accordance with the DOE Public Access Plan(http://energy.gov/downloads/doe-public-access-plan).

This work also made use of Engineering Research Center Shared Facilities supported by the Engineering Research Center Program of the National Science Foundation and the Department of Energy under NSF Award Number EEC-1041877 and the CURENT Industry Partnership Program.

$Q_{bt}^{\mathrm{U}}$, $Q_{bt}^{\mathrm{D}}$ Cost of up- and down-spinning reserve of battery $b$ during period $t$.

$P_i^{\max}, P_i^{\min}$ Maximum/minimum output of DG $i$.

$P_t^{\mathrm{W}}$, $P_t^{\mathrm{PV}}$ Wind turbine/PV power output during period $t$.

$P_{jt}$ Power consumption scheduled for demand $j$ during period $t$.

$\triangle N_t^D$ Net demand forecast error of microgrid during period $t$.

$\mu_t$, $\sigma_t$ Mean and standard deviation of $\triangle N_t^D$.

$\mathrm{PSI}^{\mathrm{req}}$ PSI requirements of microgrid operators.

$P_b^{\mathrm{C,max}}, P_b^{\mathrm{D,max}}$ Maximum charging/discharging power of battery $b$.

$SOC_{bt}^{\max}, SOC_{bt}^{\min}$ Maximum/minimum state of charge of battery $b$ during period $t$.

$\eta_b^{\mathrm{C}}, \eta_b^{D}$ Battery charging/discharging efficiency factor.

$\triangle t$ Time duration of each period.

$\tau$ Amount of time available of DGs and batteries to ramp up/down their output to deliver the reserve.

# I. Introduction

While the utilization of a microgrid for local power reliability during grid outage and emergencies is a well-known benefit, networked microgrids, defined as the aggregation of interconnected adjacent microgrids, on the other hand, offer a new, more efficient and resilient alternative to traditional individual microgrids. Due to such benefits, networked microgrids has attracted growing attention in recent years [1]-[5]. Normally, a two-layer energy management strategy for networked microgrids scheduling in distribution system has been is used. In the inner layer, each microgrid schedules its own the generation resources and loads, while the outer layer optimization coordinates the power sharing among all microgrids. From control perspective, P-Q based primary control with droop characteristics in facilitating energy transaction of the microgrids and maintaining voltage and frequency stabilities under disturbances is presented in [6].

In the existing literature, research studies on networked microgrids have been mostly focused on the optimal energy transaction strategies to meet economic objectives. However, the resiliency of microgrid and networked microgrid is rarely considered in the optimization. In fact, the most important feature of a microgrid is its ability to separate itself from the distribution utility during outage and continue supplying all or selected critical loads in its own islanded portion. Therefore, the economic benefits of networked microgrids cannot be validated without considering the system resiliency.

In view of the shortcomings of the existing networked microgrids scheduling strategies, a new scheduling strategy for both networked microgrids and independent microgrids operation considering probabilistic constraints of successful islanding is developed in this paper. Considering the uncertainty of renewable generation and power at the PCC, a new concept, probability of successful islanding (PSI), has been proposed to indicate the probability that a microgrid is maintaining adequate up- and down-spinning reserve to meet local demand and accommodate local renewable generation after instantaneously islanding from the main grid in [7]. The networked microgrids and independent microgrids are scheduled with specified PSI. The main contributions of this paper are as follows:

1) Validated the benefit of economics and resiliency of networked microgrids comparing with independent microgrids, and
2) Performed sensitivity analysis to demonstrate the impacts of correlation coefficients among the renewable generation and load of adjacent microgrids.

The rest of this paper is organized as follows. In Section II, the microgrid scheduling strategy with chance-constrained islanding capability is presented. The model is expended to networked microgrids in Section III. Case study and conclusions are given in Section IV and V.

# II. Microgrid Scheduling With Chance-Constrained Islanding Capability

## A. Component Models

The microgrid considered in this paper consists of distributed generators (e.g., diesel generators, microturbines and fuel cells), renewable generation (e.g.. wind turbines and PV panels), energy storage (e.g., battery systems) and local demands. The distributed generators are considered dispatchable units, which can be controlled by a microgrid master controller to provide both power and reserve. Depending on unit type, dispatchable units are subject to various constraints, such as, capacity limits, minimum power output limits, ramping rates, minimum on/off time, and so on. In contrast, renewable generation, such as, wind turbines and PV panels, are taken as non-dispatchable units, which depend on the meteorological conditions of wind speed, temperature and solar irradiance. Thus, renewable generation is subject to variability. Extensive research has been done on wind and PV power forecasting [8], [9]. For simplicity, we assume both wind and PV power forecast error can be modeled as independent normally distributed random variables [10]. The load forecast error is assumed to follow a normal distribution and be independent of renewable generation forecast [11]. Due to the limited size of microgrid, relatively large standard deviations are used for both renewable generation and load forecast errors.

## B. Problem Formulation

This subsection describes the model of a microgrid scheduling strategy with chance-constrained islanding capability. In the context of microgrids with dispatchable and undispatchable generation as well as electrical energy storage (e.g., batteries) integration, the objective aims at minimizing the total operation cost, including generation cost and spinning reserve cost of local resources as well as purchasing cost of energy form main grid. The objective function is shown in (1). Specifically, the first and second line is the fuel cost of DGs (including DGs start-up cost); the third line is the energy purchasing/selling cost/benefit from distribution grid; the fourth line is the battery degradation cost, the fifth and sixth lines are cost of up- and down-spinning reserve from both DGs and batteries. All terms are in mixed-integer linear form except the startup cost of generators (line 2), which can be recast into mixed-integer linear form as in [12].

$$\begin{aligned}\min \quad & \sum_{t=1}^{N_T}\sum_{i=1}^{N_G}\left[\sum_{m=1}^{N_I}\lambda_{it}(m)p_{it}(m)+A_iu_{it}\right]\\ &+\sum_{t=1}^{N_T}\sum_{i=1}^{N_G}SU_{it}\left(u_{it},u_{i,t-1}\right)\\ &+\sum_{t=1}^{N_T}\lambda_t^{\mathrm{PCC}}P_t^{\mathrm{PCC}}\\ &+\sum_{t=1}^{N_T}\sum_{b=1}^{N_B}C_{bt}\left(P_{bt}^{\mathrm{C}}+P_{bt}^{\mathrm{D}}\right)\\ &+\sum_{t=1}^{N_T}\sum_{i=1}^{N_G}\left(Q_{it}^{\mathrm{U}}R_{it}^{\mathrm{U}}+Q_{it}^{\mathrm{D}}R_{it}^{\mathrm{D}}\right)\\ &+\sum_{t=1}^{N_T}\sum_{b=1}^{N_B}\left(Q_{bt}^{\mathrm{U}}R_{bt}^{\mathrm{U}}+Q_{bt}^{\mathrm{D}}R_{bt}^{\mathrm{D}}\right) \qquad (1)\end{aligned}$$

The objective function is subject to the following constraints:

$$P_{it}=\sum_{m=1}^{N_I}p_{it}(m)+u_{it}P_i^{\min} \quad \forall i,\forall t \qquad (2)$$

$$0\le p_{it}(m)\le p_{it}^{\max}(m) \quad \forall i,\forall t,\forall m \qquad (3)$$

$$P_i^{\min}u_{it}\le P_{it}\le P_i^{\max}u_{it} \quad \forall i,\forall t \qquad (4)$$

$$R_{it}^{\mathrm{U}}\le P_i^{\max}u_{it}-P_{it} \quad \forall i,\forall t \qquad (5)$$

$$R_{it}^{\mathrm{U}}\le u_{it}R_i^{\mathrm{U,max}}\tau \quad \forall i,\forall t \qquad (6)$$

$$R_{it}^{\mathrm{D}}\le P_{it}-P_i^{\min}u_{it} \quad \forall i,\forall t \qquad (7)$$

$$R_{it}^{\mathrm{D}}\le u_{it}R_i^{\mathrm{D,max}}\tau \quad \forall i,\forall t \qquad (8)$$

$$0\le P_{bt}^{\mathrm{C}}\le P_b^{\mathrm{C,max}}u_{bt}^{\mathrm{C}} \quad \forall b,\forall t \qquad (9)$$

$$0\le P_{bt}^{\mathrm{D}}\le P_b^{\mathrm{D,max}}u_{bt}^{\mathrm{D}} \quad \forall b,\forall t \qquad (10)$$

$$u_{bt}^{\mathrm{C}}+u_{bt}^{\mathrm{D}}\le 1 \quad \forall b,\forall t \qquad (11)$$

$$SOC_{bt}=SOC_{b,t-1}+P_{bt}^{\mathrm{C}}\eta_b^{\mathrm{C}}\triangle t-P_{bt}^{\mathrm{D}}\frac{1}{\eta_b^{\mathrm{D}}}\triangle t \quad \forall b,\forall t \qquad (12)$$

$$SOC_{bt}^{\min}\le SOC_{bt}\le SOC_{bt}^{\max} \quad \forall b,\forall t \qquad (13)$$

$$P_{bt}=P_{bt}^{\mathrm{D}}-P_{bt}^{\mathrm{C}} \quad \forall b,\forall t \qquad (14)$$

$$R_{bt}^{\mathrm{U}}\le P_b^{\mathrm{D,max}}-P_{bt} \quad \forall b,\forall t \qquad (15)$$

$$R_{bt}^{\mathrm{U}}\le \eta_b^{\mathrm{D}}\left(SOC_{bt}-SOC_{bt}^{\min}\right)/\tau \quad \forall b,\forall t \qquad (16)$$

$$R_{bt}^{\mathrm{D}}\le P_b^{\mathrm{C,max}}+P_{bt} \quad \forall b,\forall t \qquad (17)$$

$$R_{bt}^{\mathrm{D}}\le 1/\eta_b^{\mathrm{C}}\left(SOC_{bt}^{\max}-SOC_{bt}\right)/\tau \quad \forall b,\forall t \qquad (18)$$

$$\begin{aligned}&\textstyle\sum_{i=1}^{N_G}P_{it}+\hat{P_t^{\mathrm{W}}}+\hat{P_t^{\mathrm{PV}}}+P_t^{\mathrm{PCC}}+\sum_{b=1}^{N_B}P_{bt}^{\mathrm{D}}\\ &\textstyle-\sum_{b=1}^{N_B}P_{bt}^{\mathrm{C}}=\sum_{j=1}^{N_D}\hat{P_{jt}} \quad \forall t \qquad (19)\end{aligned}$$

$$-\sum_{i=1}^{N_G}R_{it}^{\mathrm{D}}-\sum_{b=1}^{N_B}R_{bt}^{\mathrm{D}}\le P_t^{\mathrm{PCC}}+\triangle N_t^D\le\sum_{i=1}^{N_G}R_{it}^{\mathrm{U}}+\sum_{b=1}^{N_B}R_{bt}^{\mathrm{U}}\,\forall t \qquad (20)$$

$$\triangle N_t^D=\sum_{j=1}^{N_D}\triangle P_{jt}-\triangle P_t^{\mathrm{W}}-\triangle P_t^{\mathrm{PV}} \quad \forall t \qquad (21)$$

For DGs, constraints (2) and (3) approximate the production cost of dispatchable generators by blocks [13]. Constraint (4) forces the output of DG to be zero if it is not committed. The up-spinning reserve of DG is limited by the difference between its maximum capacity and current output in (5) and its ramping rate in (6). Similarly, the down-spinning reserve constraints are included in (7) and (8). For batteries, constraints (9) and (10) are the maximum charging/discharging power of a battery. These two states are mutually exclusive, which is ensured by (11). The battery state of charge (SOC) is defined by (12) and the limit of SOC is enforced by (13). The output power of a battery is represented in (14). Similar to DGs, the up-spinning reserve of a battery is constrained by the difference between its current SOC and minimum SOC in (15) and the difference between its maximum discharging power and current output in (16). In the same way, the down-spinning reserve constraints of a battery are included as in (17) and (18). The energy balance is enforced by (19). The spinning reserve requirement is as (20), which guarantees adequate spinning reserve for successful islanding of the microgrid considering the forecast errors of demand, wind power and PV power. The net demand forecast error $\triangle N_t^D$ is formulated in (21). Additionally, each unit or demand is subject to its own operating constraints, such as, minimum up/down time, initial condition,etc. See [14] for details about formulations of these constraints.

As mentioned in subsection II-A, we assume both wind and PV power forecast error as well as demand forecast error can be modeled as independent normally distributed random variables. Thus, the net demand forecast error $\triangle N_t^D$ also follows normal distribution, i.e., $\triangle N_t^D \sim N\left(\mu_t,\sigma_t^2\right)$. The PSI can be expressed as (22). The microgrid is considered as successfully islanded if the net demand forecast error $\triangle N_t^D \in \left[-\sum_{i=1}^{N_G}R_{it}^{\mathrm{D}}-\sum_{b=1}^{N_B}R_{bt}^{\mathrm{D}}-P_t^{\mathrm{PCC}},\ \sum_{i=1}^{N_G}R_{it}^{\mathrm{U}}+\sum_{b=1}^{N_B}R_{bt}^{\mathrm{U}} -P_t^{\mathrm{PCC}}\right]$, where $\sum_{i=1}^{N_G}R_{it}^{\mathrm{U}}+\sum_{b=1}^{N_B}R_{bt}^{\mathrm{U}}-P_t^{\mathrm{PCC}}$ stands for the redundant up-spinning reserve after islanding and $-\sum_{i=1}^{N_G}R_{it}^{\mathrm{D}}-\sum_{b=1}^{N_B}R_{bt}^{\mathrm{D}}-P_t^{\mathrm{PCC}}$ stands for the negative of the redundant down-spinning reserve after islanding. Thus, the PSI can be calculated by intergrating the probability distribution curve of $\triangle N_t^D \in \left[-\sum_{i=1}^{N_G}R_{it}^{\mathrm{D}}-\sum_{b=1}^{N_B}R_{bt}^{\mathrm{D}}-P_t^{\mathrm{PCC}},\ \sum_{i=1}^{N_G}R_{it}^{\mathrm{U}}+\sum_{b=1}^{N_B}R_{bt}^{\mathrm{U}} -P_t^{\mathrm{PCC}}\right]$, for each time interval $t$.

$$\begin{aligned}\mathrm{PSI}_t \;=\; & \mathrm{P}\Bigg(-\sum_{i=1}^{N_G}R_{it}^{\mathrm{D}}-\sum_{b=1}^{N_B}R_{bt}^{\mathrm{D}}-P_t^{\mathrm{PCC}}\le\triangle N_t^D\\ &\le\sum_{i=1}^{N_G}R_{it}^{\mathrm{U}}+\sum_{b=1}^{N_B}R_{bt}^{\mathrm{U}}-P_t^{\mathrm{PCC}}\Bigg) \qquad (22)\end{aligned}$$

The formulation of PSI considers probability distributions of forecast errors of wind, PV and loads. A multi-interval approximation of PSI is proposed in [7], which reformulates PSI into a mixed integer format. Thus, the chance-constrained programming model for microgird scheduling could be solved by mixed integer linear programming. Finally, the microgrid optimal scheduling with chance-constrained islanding capabil-ity can be formulated by substituting (20) and (21) with the linearized format of (22) and (23).

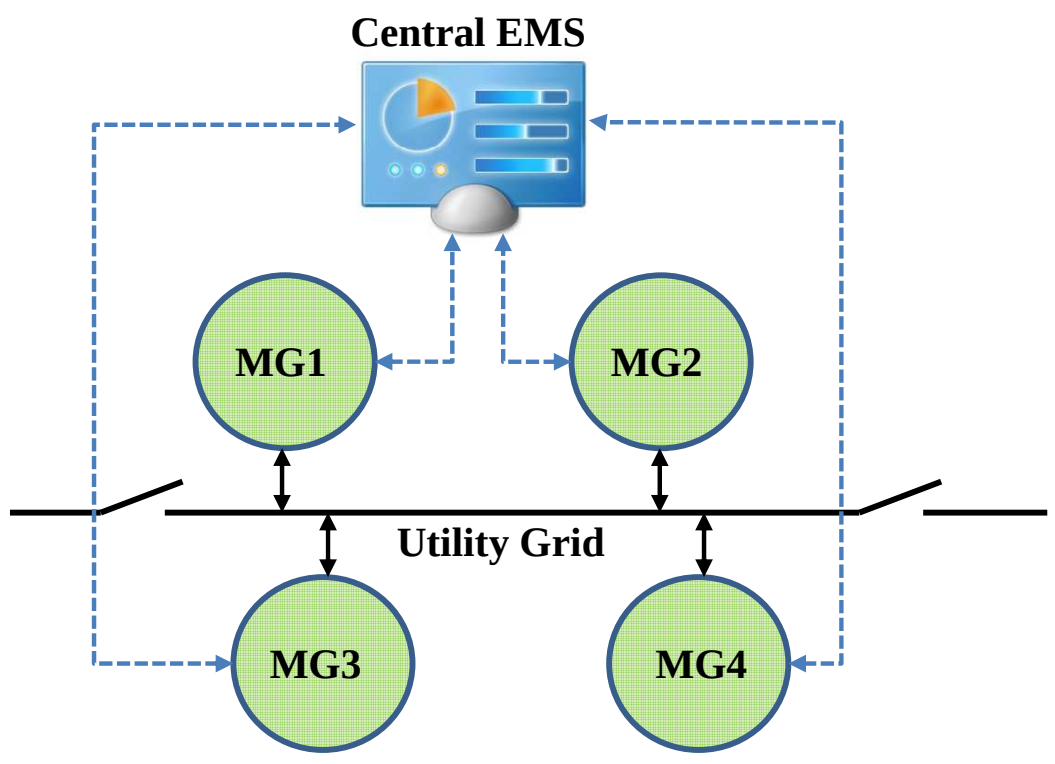

Fig. 1: Example of networked microgrids

$$\mathrm{PSI}_t \geq \mathrm{PSI}^{\mathrm{req}} \quad \forall t \tag{23}$$

The proposed chance-constrained programming model explicitly guarantees that the microgrid has adequate flexibility to meet local demand and accommodate local renewable generation after instantaneously islanding from the main grid with a certain probability specified by the microgrid operator. Thus, the resiliency of the electricity supply by the microgrid is clearly defined.

## III. Networked Microgrids Scheduling With Chance-Constrained Islanding Capability

Traditionally, each microgrid is an autonomous entity and schedules its own the generation resources and loads to maximize its own benefit. When the utility grid is faulted, each microgrid will be disconnected and performed as an autonomous island. Although multiple microgrids are physically connected, the scheduling of different microgrids is completely independent. On the other hand, interconnected adjacent microgrids can be aggregated or networked at the control and communication layer. An example of networked microgrids consisting of 4 microgrids is shown in Fig. 1. In grid-connected mode, the central EMS will schedule the 4 microgrids as a whole. When the upstream utility grid is faulted, the two switches will be opened and the 4 microgrids will formulate a single island. By networking the adjacent microgrids, better economics and resiliency are expected to be achieved comparing with independent microgrids.

In this section, we expand the resiliency-constrained scheduling model of single microgrid proposed in Section II to the case of networked microgrids. First of all, we need to substitute $P_t^{\mathrm{PCC}}$ with the summation of PCC power for all microgrids, i.e. $\sum_{n=1}^{N_M} P_{nt}^{\mathrm{PCC}}$, where $P_{nt}^{\mathrm{PCC}}$is the exchanged power at PCC for microgrid $n$ at time $t$. Secondly, we need to formulate the PSI of networked microgrids, As a precondition, the probability distribution of the net demand forecast error $\triangle N_t^D$ needs to be calculated. Same as in the previous section, we assume both wind and PV power forecast error as well as demand forecast error in a microgrid can be modeled as independent normally distributed random variables with zero mean. Due to the geographic proximity of networked microgrids, the wind power forecast errors of any two microgrids are correlated. Taking a networked microgrids consisting of 3 microgrids for example, the mean of total wind power forecast error is zero, while the deviation of total wind power forecast error can be calculated according to equation (24), where $\sigma_n^w$ is the standard deviation of wind power forecast error in microgrid $n$ and $\rho_{nn'}^w$ is the correlation coefficient between wind power forecast errors of microgrid $n$ and $n'$.

$$(\sigma^w)^2 = \begin{bmatrix} \sigma_1^w \\ \sigma_2^w \\ \sigma_3^w \end{bmatrix}^{\mathrm{T}} \begin{bmatrix} 1 & \rho_{12}^w & \rho_{13}^w \\ \rho_{21}^w & 1 & \rho_{23}^w \\ \rho_{31}^w & \rho_{32}^w & 1 \end{bmatrix} \begin{bmatrix} \sigma_1^w \\ \sigma_2^w \\ \sigma_3^w \end{bmatrix} \tag{24}$$

The standard deviation of total PV power forecast error and total demand forecast error can be calculated similarly. Since the wind and PV power forecast as well as demand forecast are independent, the total net demand forecast error of networked microgrids $\triangle N_t^D$ follows normal distribution, i.e., $\triangle N_t^D \sim N\left(\mu_t, \sigma_t^2\right)$, where $\sigma_t^2$can be easily calculated based on the result of (24). With these two modifications, the resiliency-constrained scheduling model of single microgrid has been adapted to handle the resiliency-constrained scheduling of networked microgrids.

## IV. Case Studies

In order to test the proposed networked microgrids scheduling strategy with chance-constrained islanding capability, we build a test system by connecting 3 modified ORNL Distributed Energy Control and Communication (DECC) lab microgrid on the same bus like Fig. 1. The 3 microgrids are identical. All parameters for generators, forecast wind power, PV power and demand as well as the day-ahead market prices can be found in [7]. The forecast errors of wind power and PV power are assumed to be Gaussian distribution with zero mean and 15% of standard deviation. The demand forecast error is assumed to be Gaussian distribution with zero mean and 3% of standard deviation. The analysis is conducted for a 24-hour scheduling horizon and each time interval is set to be one hour. All numerical simulations are coded in MATLAB and solved using the MILP solver CPLEX 12.2. With a prespecified duality gap of 0.1%, the running time of each case is less than 10 seconds on a 2.66 GHz Windows-based PC with 4 G bytes of RAM.

### *A. Comparing Cost of Networked Microgrids and Independent Microgrids under the Same PSI*

In order to show the benefit of networked microgrids, the total operating cost of networked microgrids and independent microgrids under the same resiliency requirements, $\mathrm{PSI}^{\mathrm{req}}$, are compared in Fig. 2. As can be seen, the operating costs of networked microgrids are always lower than that of independent microgrids under the same resiliency requirements. As the resiliency requirement $\mathrm{PSI}^{\mathrm{req}}$ increases, the economic benefit of networked microgrids becomes more significant.

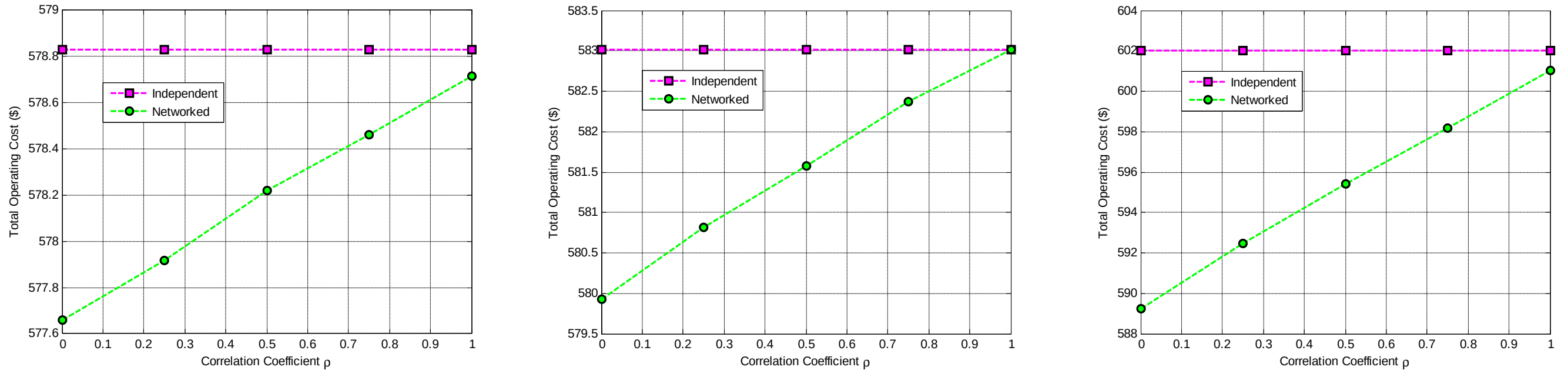

(a) $\mathrm{PSI}^{\mathrm{req}} = 0.6$ (b) $\mathrm{PSI}^{\mathrm{req}} = 0.8$ (c) $\mathrm{PSI}^{\mathrm{req}} = 0.997$

Fig. 2: Comparison of cost by networked and independent microgrids

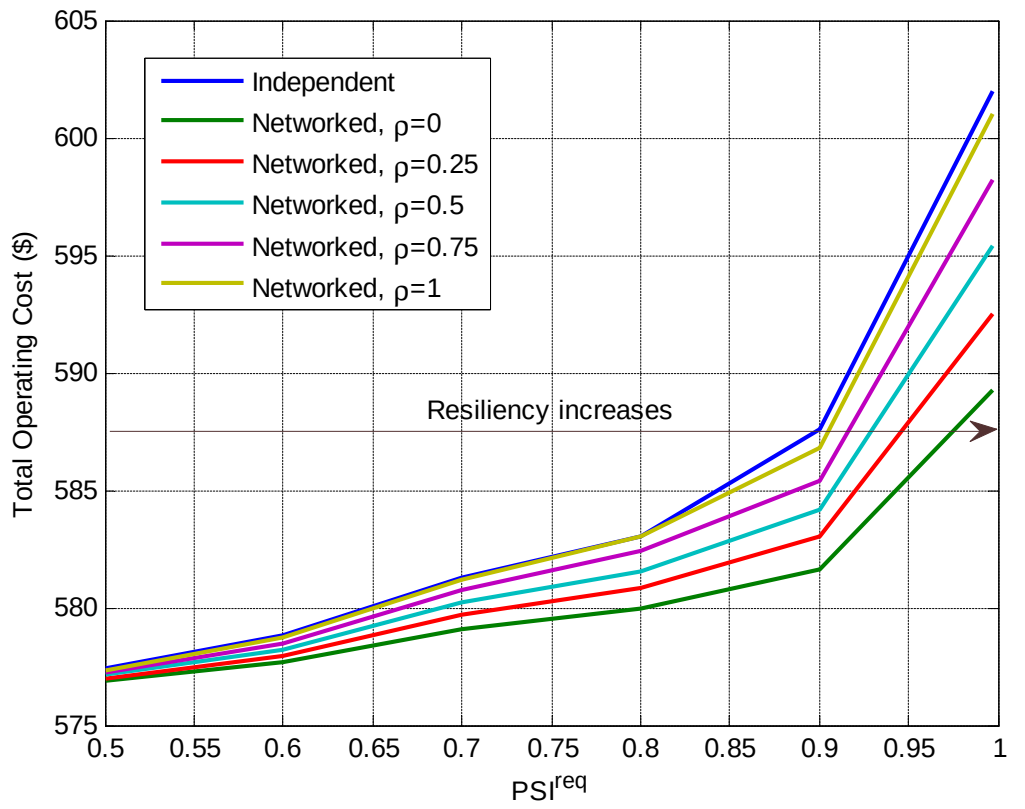

Fig. 3: Comparison of cost by networked and independent microgrids with different levels of PSI requirements

The economic benefit of networked microgrids gets smaller as the microgrids are more correlated, i.e., the correlation coefficients between different microgrids increase. Nevertheless, the economic benefit of networked microgrid is validated.

### B. Comparing PSI of Networked Microgrids and Independent Microgrids under the Same Cost

In order to show the benefit of networked microgrids in improving system resiliency, the cost by networked and independent microgrids with different levels of PSI requirements are compared in Fig. 3. As can be seen, under the same operating cost, the networked microgrids always have higher resiliency. This effect is much more obvious when the microgrids are less correlated. This clearly validates the resiliency benefit of networked microgrids.

## V. Conclusions

In this paper, we modified the resiliency-constrained scheduling model of single microgrid to handle the case of networked microgrids. The model explicitly defines the resiliency of a microgrid and networked microgrids considering islanding situations and forecast uncertainties. Numerical simulations validated the benefit of networked microgrids in the aspect of economics and resiliency comparing with independent microgrids. In addition, the impact of correlation coefficients among the renewable generation and load of adjacent microgrids has been studied as well.

## References

[1] Z. Wang, B. Chen, J. Wang, M. M. Begovic and C. Chen, "Coordinated Energy Management of Networked Microgrids in Distribution Systems," *IEEE Trans. Smart Grid*, vol. 6, no. 1, pp. 45-53, Jan. 2015.

[2] Z. Wang, B. Chen, J. Wang and C. Chen, "Networked Microgrids for Self-Healing Power Systems," *IEEE Trans. Smart Grid*, vol. 7, no. 1, pp. 310-319, Jan. 2016.

[3] A. Arif and Z. Wang, "Networked microgrids for service restoration in resilient distribution systems," *IET Gener. Transm. Distrib.*, vol. 11, no. 14, pp. 3612-3619, 9 28 2017.

[4] A. Hussain, V. H. Bui and H. M. Kim, "A Resilient and Privacy-Preserving Energy Management Strategy for Networked Microgrids," *IEEE Trans. Smart Grid*, vol. . PP, no. 99, pp. 1-1.

[5] A. Parisio, C. Wiezorek, T. Kyntäjä, J. Elo, K. Strunz and K. H. Johansson, "Cooperative MPC-Based Energy Management for Networked Microgrids," *IEEE Trans. Smart Grid*, vol. 8, no. 6, pp. 3066-3074, Nov. 2017.

[6] F. Zhang, H. Zhao and M. Hong, "Operation of networked microgrids in a distribution system," in CSEE Journal of Power and Energy Systems, vol. 1, no. 4, pp. 12-21, Dec. 2015.

[7] G. Liu, M. Starke, B. Xiao, X. Zhang, K. Tomsovic, "Microgrid optimal scheduling with chance-constrained islanding capability," *Electr. Power Syst. Res.*, vol. 145, pp. 197-206, Apr. 2017.

[8] J. Carta, P. Ramrez, and S. Velázquez, " A review of wind speed probability distributions used in wind energy analysis," *Renew. Sustain. Energy Rev.*, vol. 13, no. 5, pp. 933-955, Jun. 2009.

[9] P. Bacher, H. Madsen, and H. A. Nielsen, "Online short-term solar power forecasting," *Solar Energy*, vol. 83, pp. 1772-1783, 2009.

[10] B. Zhao, Y. Shi, X. Dong, W. Luan, and J. Bornemann, "Short-Term Operation Scheduling in Renewable-Powered Microgrids: A Duality-Based Approach," *IEEE Trans. Sustain. Energy*, vol. 5, no. 1, pp. 209-217, Jan. 2014.

[11] M. Q. Wang, and H. B. Gooi, "Spinning Reserve Estimation in Microgrids," *IEEE Trans. Power Syst.*,, vol.26, no.3, pp.1164-1174, Aug. 2011.

[12] M. Ortega-Vazquez, "Optimizing the spinning reserve requirements," in *Sch. Elect. Electron. Eng.*. Manchester, U.K.: Univ. Manchester, 2006, pp. 1-219 [Online]. Available: http://www.eee.manchester.ac.uk/research/groups/eeps/publications/reportstheses/aoe/ortega-vazquez_PhD_2006.pdf.

[13] F. Aminifar, M. Fotuhi-Firuzabad, and M. Shahidehpour, "Unit commitment with probabilistic spinning reserve and interruptible load considerations," *IEEE Trans. Power Syst.*, vol. 24, no. 1, pp. 388-397, Feb. 2009.

[14] M. Carrión, and J. M. Arroyo, "A computationally efficient mixed-integer linear formulation for the thermal unit commitment problem," *IEEE Trans. Power Syst.*, vol. 21, no. 3, pp. 1371-1378, Aug. 2006.